\documentclass[12pt]{article}
\usepackage{graphicx}
\usepackage{amsmath,amssymb,amsthm,amsfonts}
\newtheorem{thm}{Theorem}[section]

\newcommand{\R}{{\mathbb{R}}}
\newcommand{\Z}{{\mathbb{Z}}}

\newcommand{\1}{\partial}
\newcommand{\2}{\overline}
\newcommand{\3}{\varepsilon}
\newcommand{\4}{\widetilde}

\def\ni{\noindent}

\begin{document}
\title{A pseudolocality theorem for Ricci flow}
\author{Shu-Yu Hsu\\
Department of Mathematics\\
National Chung Cheng University\\
168 University Road, Min-Hsiung\\
Chia-Yi 621, Taiwan, R.O.C.\\
e-mail: syhsu@math.ccu.edu.tw}
\date{Oct 6, 2010}
\smallbreak \maketitle
\begin{abstract}
In this paper we will give a simple proof of a modification of a result on 
pseudolocality for the Ricci flow by P.~Lu \cite{L} without 
using the pseudolocality theorem 10.1 of Perelman \cite{P1}. We also
obtain an extension of a result of Hamilton \cite{H3}
on the compactness of a sequence of complete pointed Riemannian 
manifolds $\{(M_k,g_k(t),x_k)\}_{k=1}^{\infty}$ evolving 
under Ricci flow with uniform bounded sectional curvatures on $[0,T]$ and 
uniform positive lower bound on the injectivity radii at $x_k$ with respect 
to the metric $g_k(0)$.
\end{abstract}

\vskip 0.2truein

Key words: pseudolocality, complete Riemannian manifold, Ricci flow, 
locally bounded Riemannian curvature

AMS Mathematics Subject Classification: Primary 58J35, 53C44 Secondary 35K55
\vskip 0.2truein
\setcounter{equation}{0}
\setcounter{section}{-1}

\setcounter{equation}{0}
\setcounter{thm}{0}

A time dependent metric $g_{ij}(t)$ on an $n$-dimensional manifold $M$ 
is said to evolve by the Ricci flow on $(0,T)$ if it satisfies
\begin{equation}
\frac{\1}{\1 t}g_{ij}=-2R_{ij}
\end{equation}
on $(0,T)$ where $R_{ij}(x,t)$ is the Ricci tensor with respect to
the metric $g_{ij}(t)$. In 1982 R.S.~Hamilton \cite{H1} used Ricci flow 
to prove that any compact $3$-dimensional Riemannian  manifold
with strictly positive Ricci curvature also admits a metric of constant 
positive curvature. Recently there are many research on Ricci flow by
A.~Chau, L.F.~Tam and C.~Yu \cite{CTY}, P.~Daskalopoulos, R.S.~Hamilton 
and  N.~Sesum \cite{DHS}, \cite{H1}, \cite{H3}, S.Y.~Hsu 
\cite{Hs1}, \cite{Hs2}, 
\cite{Hs3}, B.~Kleiner and J.~Lott \cite{KL}, J.~Morgan and G.~Tang \cite{MT}, 
L.~Ni \cite{N}, G.~Perelman \cite{P1}, \cite{P2}, S.~Kuang and 
Q.S.~Zhang \cite{KZ}, 
\cite{Z} etc. Interested readers can read the survey article \cite{H2} by 
R.S.~Hamilton and the book Hamilton's Ricci flow \cite{CLN} 
by B.~Chow, P.~Lu, and L.~Ni for more results on Ricci flow.

In a recent paper \cite{L} P.~Lu proved the following pseudolocality 
theorem for Ricci flow.

\begin{thm}
For any $n\in\Z^+$ and $\delta>0$ there exists a constant $\3_0>0$ with the
following property. For any $r_0>0$ and $0<\3<\3_0$ suppose $(M,g(t))$ 
is an $n$-dimensional complete solution to the Ricci flow on $[0,(\3 r_0)^2]$ 
with bounded sectional curvature, and assume that there exists $x_0\in M$ 
such that 
\begin{equation}
|Rm|(x,0)\le r_0^{-2}\quad\forall x\in B_{g(0)}(x_0,r_0)
\end{equation}
and
\begin{equation}
\mbox{Vol}_{g(0)}(B_{g(0)}(x_0,r_0))\ge\delta r_0^n.
\end{equation}
Then 
\begin{equation}
|Rm|(x,t)\le (\3 r_0)^{-2}\quad\forall x\in B_{g(t)}(x_0,\3 r_0),
0\le t\le (\3 r_0)^2.
\end{equation}
\end{thm}

As observed by P.~Lu \cite{L} Theorem 1 is implied by the following 
theorem.

\begin{thm}
For any $n\in\Z^+$ and $\delta>0$ there exists a constant $\3_0>0$ with the
following property. For any $r_0>0$ and $0<\3<\3_0$ suppose $(M,g(t))$ 
is an $n$-dimensional complete solution to the Ricci flow on $[0,(\3 r_0)^2]$ 
with bounded sectional curvature, and assume that there exists $x_0\in M$ 
such that (2) and (3) hold. Then 
\begin{equation}
|Rm|(x,t)\le (\3 r_0)^{-2}\quad\forall x\in B_{g(0)}(x_0,e^{n-1}\3 r_0),
0\le t\le (\3 r_0)^2.
\end{equation}
\end{thm}

The proof of Theorem 2 in \cite{L} uses the pseudolocality theorem Theorem 
10.1 of Perelman \cite{P1}. However a careful examination of the proof
of Theorem 10.1 of Perelman \cite{P1} shows that the proof of Theorem 10.1 
of \cite{P1} is not correct. The reason is as follows. In the proof of 
Theorem 10.1 of \cite{P1} Perelman constructed a sequence of pointed
Ricci flow $(M_k, g_k(t), (x_{0,k},0))$, $0\le t\le\3_k$, with 
$\3_k\to 0$ as $k\to\infty$ and a sequence $\delta_k\to 0$ as 
$k\to\infty$ that satisfies
\begin{equation}
|Rm_{g_k}|(x,t)\le\alpha t^{-1}+2\3_k^{-2}\quad\forall d_{g_k(t)}
(x,x_{0,k})\le\3_k,0\le t\le\3_k^2
\end{equation}
for some constant $\alpha>0$ and a sequence $(x_k,t_k)$ with
$0<t_k\le\3_k^2$ and $d_{g_k(t_k)}(x_{0,k},x_k)<\3_k$ such that
$$
|Rm_{g_k}(x_k,t_k)|>\alpha t_k^{-1}+\3_k^{-2}.
$$
Perelman \cite{P1} also constructed a sequence $(\2{x}_k,\2{t}_k)$ with
$$
0<\2{t}_k\le\3_k^2\quad\mbox{  and }\quad d_{g_k(\2{t}_k)}(x_{0,k},\2{x}_k)
<(2A_k+1)\3_k
$$ 
where $A_k=1/(100n\3_k)$ such that
$$
|Rm_{g_k}(x,t)|\le 4Q_k\quad\forall d_{g_k(\2{t}_k)}(x,\2{x}_k)
<\frac{1}{10}Q_k^{-\frac{1}{2}},\2{t}_k-\frac{1}{2}\alpha Q_k^{-1}
\le t\le\2{t}_k
$$
where $Q_k=|Rm_{g_k(t)}(\2{x}_k,\2{t}_k)|$. On the third paragraph on 
P.26 of \cite{P1} Perelman claimed that the sequence of metrics
$\hat{g}_k(t)=\frac{1}{2\2{t}_k}g(2\2{t}_kt)$ converges 
to some solution of Ricci flow $\hat{g}_{\infty}(t)$ on $0\le t\le 1/2$ 
as $k\to\infty$. Perelman then concluded that there is a contradiction
to the logarithmic Sobolev inequality on $\R^n$ by passing to the limit 
a rescaled version of the equation on P.26 of \cite{P1} for $t=0$ as 
$k\to\infty$. However by (6),
\begin{equation*}
|Rm_{\hat{g}_k}|(x,t)\le 2\2{t}_k|Rm_{g_k}|(x,2\2{t}_kt)
\le 2\2{t}_k(\alpha (2\2{t}_kt)^{-1}+2\3_k^{-2})=\alpha t^{-1}+4
\end{equation*}
for any $d_{\hat{g}_k(t)}(x,x_0)\le\3_k/\sqrt{2\2{t}_k}$ and 
$0\le t\le 1/2$. Hence $|Rm_{\hat{g}_k}|(x,t)$ are not uniformly 
bounded near $t=0$. Thus one cannot apply Hamilton's compactness
theorem \cite{H3} to conclude that the sequence $\hat{g}_k(t)$
converges to some solution of Ricci flow $\hat{g}_{\infty}(t)$ on 
$0\le t\le 1/2$ as $k\to\infty$. It is also not known why one can 
pass to the limit for the inequality on P.26 of \cite{P1} as 
$k\to\infty$. 

Hence the proof of Theorem 10.1 of \cite{P1} is 
not correct and the validity of Theorem 10.1 of \cite{P1} is not 
known. On the other hand in the more detailed explanation of the 
proof of Theorem 10.1 of \cite{P1} on P.179 of \cite{CCG} it is  
hard to check that the function $\psi_i^2(x)
=(2\pi)^{-n/2}e^{-\4{f}_i(x,0)}$ defined there
belong to $W^{1,2}$ which is the required condition for the validity
of the logarithmic Sobolev inequality (Theorem 22.16 of \cite{CCG})
for manifolds satisfying the isoperimetric inequality. 

In this paper we prove that under a mild additional hypothesis 
Theorem 2 holds without using Theorem 10.1 of \cite{P1}. More 
specifically we will prove that the following result holds.

\begin{thm}
For any $n\in\Z^+$, $C_0>0$ and $\delta>0$ there exists a constant 
$\3_0>0$ with the following property. For any $r_0>0$ and $0<\3<\3_0$ 
suppose $(M,g(t))$ is an $n$-dimensional complete solution to the Ricci 
flow on $[0,(\3 r_0)^2]$ with bounded sectional curvature, and assume 
that there exists $x_0\in M$ such that (2), (3) and
\begin{equation}
|Rm|(x,t)\le\frac{C_0}{t}\quad\forall x\in B_{g(0)}(x_0,r_0),
0<t\le (\3r_0)^2
\end{equation}
hold. Then (4) holds.
\end{thm}

Similar to \cite{L} Theorem 3 is implied by the following theorem.

\begin{thm}
For any $n\in\Z^+$, $C_0>0$ and $\delta>0$ there exists a constant 
$\3_0>0$ with the following property. For any $r_0>0$ and $0<\3<\3_0$ 
suppose $(M,g(t))$ is an $n$-dimensional complete solution to the Ricci 
flow on $[0,(\3 r_0)^2]$ with bounded sectional curvature, and assume 
that there exists $x_0\in M$ such that (2), (3) and (7) hold. Then 
(5) holds.
\end{thm}
 
\vskip 0.1truein
{\ni{\it Proof of Theorem 4}:}
By rescaling the metric by $1/r_0^2$ we may assume without loss of 
generality that $r_0=1$. Suppose the theorem is not true. Then there 
exist $n\in\Z^+$, $C_0>0$, $\delta>0$, a sequence of positive numbers 
$0<\3_k<e^{1-n}/5$ with $\3_k\to 0$ as $k\to\infty$, 
and a sequence of $n$-dimensional complete manifolds $(M_k,g_k(t))$, 
$0\le t\le\3_k^2$, with $g_k$ satisfying the Ricci flow on 
$[0,\3_k^2]$ with bounded sectional curvature, a sequence $x_{0,k}
\in M_k$ and a sequence $(x_k,t_k)\in B_{g_k(0)}(x_{0,k},e^{n-1}\3_k)\times 
(0,\3_k^2]$ such that 
\begin{equation}
|Rm_{g_k(0)}|(x,0)\le 1\quad\forall x\in B_{g_k(0)}(x_{0,k},1),
\end{equation}
\begin{equation}
|Rm_{g_k(t)}|(x,t)\le\frac{C_0}{t}\quad\forall x\in B_{g_k(0)}
(x_{0,k},1),0<t\le\3_k^2,
\end{equation}
\begin{equation}
\mbox{Vol}_{g_k(0)}(B_{g_k(0)}(x_{0,k},1))\ge\delta
\end{equation}
and
\begin{equation}
|Rm_{g_k}|(x_k,t_k)>\3_k^{-2}
\end{equation}
holds for all $k\in\Z^+$. 
\vspace{6pt}

Then by \cite{L} we have the following result.
\vspace{6pt}

\noindent {\bf Claim 1} (Claim A of \cite{L}): For any $k\in\Z^+$, 
there exists $(\2{x}_k,\2{t}_k)\in B_{g_k(0)}(x_{0,k},(2A_k+e^{n-1})\3_k)
\times (0,\3_k^2]$ with $Q_k=|Rm|(\2{x}_k,\2{t}_k)>\3_k^{-2}$ such that 
\begin{equation}
|Rm|(x,t)\le 4Q_k\quad\forall (x,t)\in B_{g_k(0)}
(\2{x}_k,A_kQ_k^{-\frac{1}{2}})\times (0,\2{t}_k]
\end{equation}
where $A_k=1/(100n\3_k)$. 

\vspace{6pt}

Let $\hat{g}_k(t)=Q_kg_k(t/Q_k)$ and
$\hat{t}_k=\2{t}_kQ_k$. By passing to a subsequence if 
necessary we may assume without loss of generality that
$d_{g_k(0)}(\2{x}_k,x_{0,k})<1/10$ and $A_k\ge 2$ for all 
$k\in\Z^+$ and
$$
T=\lim_{k\to\infty}\hat{t}_k\in [0,C_0]
$$
exists. Then 
\begin{equation}\left\{\begin{aligned}
&|Rm_{\hat{g}_k}|(\2{x}_k,\hat{t}_k)=1\\
&|Rm_{\hat{g}_k}|(x,t)\le 4\qquad\,\,\forall (x,t)\in B_{\hat{g}_k(0)}
(\2{x}_k,A_k)\times (0,\hat{t}_k]\\
&|Rm_{\hat{g}_k}|(x,0)\le Q_k^{-1}\quad\forall x\in B_{\hat{g}_k(0)}
(\2{x}_k,Q_k^{1/2}/2)
\end{aligned}\right.
\end{equation}
hold for all $k\in\Z^+$.
By (8), (10), and the Bishop volume comparison theorem there exists a 
constant $\delta_1>0$ such that
\begin{equation*}
\mbox{Vol}_{g_k(0)}(B_{g_k(0)}(x_{0,k},1/10))\ge\delta_1
\quad\forall k\in\Z^+.
\end{equation*}
Since  $B_{g_k(0)}(x_{0,k},1/10)\subset B_{g_k(0)}(\2{x}_k,1/4)$, 
\begin{equation}
\mbox{Vol}_{g_k(0)}(B_{g_k(0)}(\2{x}_k,1/4))\ge\delta_1\quad\forall 
k\in\Z^+.
\end{equation}
Since  $B_{g_k(0)}(\2{x}_k,1/4)\subset B_{g_k(0)}(x_{0,k},1)$,
by (8), (14) and Lemma 1 of \cite{L} (cf.  \cite{W} and Theorem 4.10 
of \cite{H3}) there exist constants $\delta_0>0$ and $0<r_0<1/4$ such that 
\begin{equation*}
\mbox{Area}_{g_k(0)}(\1\Omega)^n\ge (1-\delta_0)
\mbox{Vol}_{g_k(0)}(\Omega)^{n-1}
\end{equation*} 
holds for any regular domain $\Omega\subset B_{g_k(0)}(\2{x}_k,r_0)$ and 
$k\in\Z^+$. Hence 
\begin{equation}
\mbox{Area}_{\hat{g}_k(0)}(\1\Omega)^n\ge (1-\delta_0)
\mbox{Vol}_{\hat{g}_k(0)}(\Omega)^{n-1}
\end{equation} 
holds for any regular domain $\Omega\subset B_{\hat{g}_k(0)}(\2{x}_k,Q_k^{1/2}r_0)$ 
and $k\in\Z^+$. Since $Q_k\to\infty$ as $k\to\infty$, without loss of 
generality we may assume that $Q_k^{1/2}r_0>1$ for all $k\in\Z^+$. Then by (15) 
there exists a positive constant $\delta_2>0$ such that
\begin{equation}
\mbox{Vol}_{\hat{g}_k(0)}(B_{\hat{g}_k(0)}(\2{x}_k,1))\ge\delta_2
\quad\forall k\in\Z^+.
\end{equation}
We now divide the proof into two cases.

\vspace{6pt}

\noindent $\underline{\text{\bf Case 1}}$: $T=0$

\vspace{6pt}
This case can be shown to be impossible by the same argument as the proof
of case 3 on P.8--9 of \cite{L} using Theorem 8.3 of
\cite{P1} and a modification of the argument of Perelman \cite{P1}. 
For the sake of completeness we will give a simple different
proof here.
For any $k\in\Z^+$ let $\delta_{\2{x}_k}$ be the delta mass at $\2{x}_k$
and $\eta_k$ be the solution of
\begin{equation}\left\{\begin{aligned}
&\eta_{k,t}+\Delta_{\hat{g}_k(t)}\eta_k+C_1\eta_k=0
\quad\mbox{ in }B_{\hat{g}_k(0)}(\2{x}_k,1)\times (0,\hat{t}_k]\\
&\eta_k(x,\hat{t}_k)=\delta_{\2{x}_k}
\qquad\qquad\qquad\mbox{ in }B_{\hat{g}_k(0)}(\2{x}_k,1)\\
&\eta_k(x,t)=0\qquad\qquad\qquad\quad\,\mbox{ on }\1 B_{\hat{g}_k(0)}
(\2{x}_k,1)\times (0,\hat{t}_k]
\end{aligned}\right.
\end{equation}
where $C_1=64+4n(n-1)$. Then by the maximum principle $\eta_k\ge 0$ in
$B_{\hat{g}_k(0)}(\2{x}_k,1)\times (0,\hat{t}_k]$ and $\1\eta_k/\1\nu\ge 0$
on $\1 B_{\hat{g}_k(0)}(\2{x}_k,1)\times (0,\hat{t}_k]$ where
$\1/\1\nu$ is the derivative with respect to the unit inward normal
$\nu$ on $\1 B_{\hat{g}_k(0)}(\2{x}_k,1)\times (0,\hat{t}_k]$. We extend
$\eta_k$ by letting $\eta_k =0$ on $(\2{B}_{\hat{g}_k(0)}(\2{x}_k,1)\times
[\hat{t}_k,\infty))\setminus\{(\2{x}_k,\hat{t}_k)\}$ and we extend 
$\hat{g}_k(t)$ by letting $\hat{g}_k(t)=\hat{g}_k(\hat{t}_k)$ for all
$t\ge\hat{t}_k$. Then $\eta_k$ satisfies
$$
\eta_{k,t}+\Delta_{\hat{g}_k(t)}\eta_k+C_1\eta_k=0\quad\mbox{ in }
(\2{B}_{\hat{g}_k(0)}(\2{x}_k,1)\times
(0,\infty))\setminus\{(\2{x}_k,\hat{t}_k)\}.
$$
Now by (1), (13), and (17),
\begin{align}
\frac{\1}{\1 t}\biggl (\int_{B_{\hat{g}_k(0)}(\2{x}_k,1)}\eta_k
\,d\hat{V}_k(t)\biggr)
=&\int_{B_{\hat{g}_k(0)}(\2{x}_k,1)}(\eta_{k,t}-R_{\hat{g}_k}\eta_k)
\,d\hat{V}_k(t)\nonumber\\
=&\int_{B_{\hat{g}_k(0)}(\2{x}_k,1)}(-\Delta_{\hat{g}_k(t)}\eta_k
-(C_1+R_{\hat{g}_k})\eta_k)\,d\hat{V}_k(t)\nonumber\\
\ge&-C_4\int_{B_{\hat{g}_k(0)}(\2{x}_k,1)}\eta_k\,d\hat{V}_k(t)
+\int_{\1 B_{\hat{g}_k(0)}(\2{x}_k,1)}\frac{\1\eta_k}{\1\nu}
\,d\sigma_k(t)\nonumber\\
\ge&-C_4\int_{B_{\hat{g}_k(0)}(\2{x}_k,1)}\eta_k\,d\hat{V}_k(t)
\end{align}
for any $0\le t<\hat{t}_k$ where $C_4=64+8n(n-1)$. Integrating (18)
over $(0,\hat{t}_k)$,
\begin{equation}
\int_{B_{\hat{g}_k(0)}(\2{x}_k,1)}\eta_k (x,t)\,d\hat{V}_k(t)
\le e^{C_4\hat{t}_k}\quad\forall 0\le t<\hat{t}_k.
\end{equation}
Hence by (13), (19) and the parabolic Schauder estimates \cite{LSU} 
(cf.\cite{CTY},\cite{KZ}) there exists a constant $C_2>0$ such that
\begin{equation}
|\1\eta_k/\1\nu|\le C_2\quad\mbox{ on }\1 B_{\hat{g}_k(0)}
(\2{x}_k,1)\times [0,\hat{t}_k]\quad\forall k\in\Z^+.
\end{equation}   
Since the curvature $Rm_{\hat{g}_k}$ satisfies (cf. \cite{L}, \cite{H2}),
\begin{equation*}
(|Rm_{\hat{g}_k}|^2)_t\le\Delta_{\hat{g}_k(t)}|Rm_{\hat{g}_k}|^2
-2|\nabla_{\hat{g}_k}Rm_{\hat{g}_k}|^2+16|Rm_{\hat{g}_k}|^3,
\end{equation*} 
by (1), (13), (17), and (20) we have
\begin{align}
&\frac{\1}{\1 t}\biggl(\int_{B_{\hat{g}_k(0)}(\2{x}_k,1)}
|Rm_{\hat{g}_k}|^2\eta_k\,d\hat{V}_k(t)\biggr)\nonumber\\
=&\int_{B_{\hat{g}_k(0)}(\2{x}_k,1)}[(|Rm_{\hat{g}_k}|^2)_t\eta_k
+|Rm_{\hat{g}_k}|^2\eta_{k,t}-R_{\hat{g}_k}|Rm_{\hat{g}_k}|^2\eta_k]
\,d\hat{V}_k(t)\nonumber\\
\le&\int_{B_{\hat{g}_k(0)}(\2{x}_k,1)}[(\Delta_{\hat{g}_k(t)}
|Rm_{\hat{g}_k}|^2+64|Rm_{\hat{g}_k}|^2)\eta_k
+|Rm_{\hat{g}_k}|^2\eta_{k,t}-R_{\hat{g}_k}|Rm_{\hat{g}_k}|^2\eta_k]
\,d\hat{V}_k(t)\nonumber\\
\le&\int_{B_{\hat{g}_k(0)}(\2{x}_k,1)}|Rm_{\hat{g}_k}|^2
(\eta_{k,t}+\Delta_{\hat{g}_k(t)}\eta_k+C_1\eta_k)\,d\hat{V}_k(t)
+16\int_{\1 B_{\hat{g}_k(0)}(\2{x}_k,1)}
\frac{\1\eta_k}{\1\nu}\,d\sigma_k(t)\nonumber\\
\le&16C_2|\1 B_{\hat{g}_k(0)}(\2{x}_k,1)|
\end{align}
for any $0\le t<\hat{t}_k$ where $R_{\hat{g}_k}$, $d\hat{V}_k(t)$, 
$d\sigma_k(t)$ are the scalar curvature, volume element, and surface 
element with respect to the metric $\hat{g}_k(t)$. 

By (13), (16), and 
Cheeger-Gromov's compactness theorem (\cite{G}, \cite{Pe}) the sequence 
of pointed manifold $(M_k,\hat{g}_k(0),\2{x}_k)$ has a subsequence
which we may assume without loss of generality to be the sequence itself
that converges to some pointed manifold $(M_0,\hat{g}_0,\2{x}_0)$ as 
$k\to\infty$ (cf. \cite{Pe},\cite{H3}). Then there exists a constant 
$C_3>0$ such that
\begin{equation}
|\1 B_{\hat{g}_k(0)}(\2{x}_k,1)|\le\frac{C_3}{16C_2}\quad\forall
k\in\Z^+.
\end{equation}
Hence by (13), (17), (21) and (22),
\begin{align}
1=&\lim_{t\nearrow\hat{t}_k}\int_{B_{\hat{g}_k(0)}(\2{x}_k,1)}
|Rm_{\hat{g}_k}(x,t)|^2\eta_k(x,t)\,d\hat{V}_k(t)\nonumber\\
\le&\int_{B_{\hat{g}_k(0)}(\2{x}_k,1)}|Rm_{\hat{g}_k}(x,0)|^2\eta_k(x,0)
\,d\hat{V}_k(0)+C_3\hat{t}_k\nonumber\\
\le&Q_k^{-2}\int_{B_{\hat{g}_k(0)}(\2{x}_k,1)}\eta_k(x,0)\,d\hat{V}_k(0)
+C_3\hat{t}_k
\end{align}
By (19) and (23),
\begin{equation}
1\le Q_k^{-2}e^{C_4\hat{t}_k}+C_3\hat{t}_k\quad\forall k\in\Z^+.
\end{equation}
Letting $k\to\infty$ in (24) we get $1\le 0$ and contradiction arises.
Hence $T>0$.

\vspace{6pt}

\noindent $\underline{\text{\bf Case 2}}$: $T\in (0,C_0]$

\vspace{6pt}
By passing to a subsequence if necessary we may assume without loss 
of generality that
$T/2<\hat{t}_k<3T/2$ and $A_k>e^{9nT}$ for all $k\in\Z^+$. By (1) and (13),
\begin{equation}
e^{-6nT}\hat{g}_k(0)\le \hat{g}_k(t)\le e^{6nT}\hat{g}_k(0)\quad
\mbox{ in }B_{\hat{g}_k(0)}(\2{x}_k, A_k)\times (0,\hat{t}_k]
\quad\forall k\in\Z^+.
\end{equation}
Then 
\begin{equation}
B_{\hat{g}_k(0)}(\2{x}_k,1)\subset B_{\hat{g}_k(t)}(\2{x}_k,e^{3nT})
\subset B_{\hat{g}_k(0)}(\2{x}_k,e^{9nT})\quad\forall 0\le t\le \hat{t}_k,
k\in\Z^+.
\end{equation} 
Hence by (16) and (26),
\begin{equation}
Vol_{\hat{g}_k(0)}(B_{\hat{g}_k(T/2)}(\2{x}_k,e^{3nT}))\ge\delta_2\quad\forall
k\in\Z^+.
\end{equation}
Now by (1) and (13),
\begin{equation}
\left|\frac{\1}{\1 t}(\log (d\hat{V}_k(t)))\right|\le 4n(n-1)\quad
\mbox{ in }B_{\hat{g}_k(0)}(\2{x}_k,A_k)\times (0,\hat{t}_k]
\quad\forall k\in\Z^+.
\end{equation}
Hence by (26) and (28),
\begin{equation}
Vol_{\hat{g}_k(T/2)}(B_{\hat{g}_k(T/2)}(\2{x}_k,e^{3nT}))\ge
e^{-2n(n-1)T}Vol_{\hat{g}_k(0)}(B_{\hat{g}_k(T/2)}(\2{x}_k,e^{3nT}))
\quad\forall k\in\Z^+.
\end{equation}
By (27) and (29),
\begin{equation}
Vol_{\hat{g}_k(T/2)}(B_{\hat{g}_k(T/2)}(\2{x}_k,e^{3nT}))\ge
e^{-2n(n-1)T}\delta_2\quad\forall k\in\Z^+.
\end{equation}
Let $M_k=B_{\hat{g}_k(0)}(\2{x}_k,A_k)$.
By (13) and (30) the injectivity radii of $(M_k,\hat{g}_k(T/2))$ 
at $\2{x}_k$ are uniformly bounded below by some positive constant for all
$k\in\Z^+$. Hence by (13) and the Hamilton compactness theorem \cite{H3}
there exists a subsequence of $(M_k,\hat{g}_k(t),(\2{x}_k,T/2))$ which we 
may assume without loss of generality to be the sequence itself that 
converges to some pointed complete manifold 
$(M_{\infty},g_{\infty}(t),(x_{\infty},T/2))$, $0<t\le T$, as 
$k\to\infty$. $g_{\infty}$ satisfies the Ricci flow equation (1) with 
\begin{equation}
|Rm_{\infty}|(x,t)\le 4\quad\forall x\in M_{\infty}, 0<t\le T
\end{equation}
where $Rm_{\infty}(x,t)$ is the Riemmannian curvature of $M_{\infty}$
with respect to the metric $g_{\infty}(t)$. Let
$$
h(t)=\sup_{x\in M_{\infty}}|Rm_{\infty}|(x,t).
$$
Then $0\le h(t)\le 4$ on $(0,T]$.
Since $|Rm_{\hat{g}_k(\hat{t}_k)}(\2{x}_k,\hat{t}_k)|=1$ for all $k\in\Z^+$, 
\begin{equation}
|Rm_{\infty}(x_{\infty},T)|=1. 
\end{equation}
Hence $h(T)\ge 1$. By continuity
there exists a constant $\delta_1>0$ such that $h(t)>0$ on 
$[T-\delta_1,T]$. Let $T_1\ge 0$ be the minimal time such that 
$h(t)>0$ for any $t\in (T_1,T]$. Then $T_1<T-\delta_1$.
Suppose $T_1>0$. Then
\begin{equation}
Rm_{\infty}(x,T_1)\equiv 0\quad\mbox{ on }M_{\infty}.
\end{equation}
By (1) and (33),
\begin{equation}
Rm_{\infty}(x,t)\equiv 0\quad\mbox{ on }M_{\infty}\times (0,T].
\end{equation}
By (32) and (34) contradiction arises. Hence $T_1=0$.
Let $\delta_2=\inf_{0<t\le T}h(t)$.

We now divide the proof of case 2 into two subcases:

\noindent $\underline{\mbox{\bf Case (a)}}$: $\delta_2=0$.

Then there exists a sequence $\{s_i\}_{i=1}^{\infty}$ with $0<s_i<T$ for 
all $i\in\Z^+$ and $s_i\to 0$ as $i\to\infty$ such that $h(s_i)\to 0$
as $i\to\infty$. Since the curvature $Rm_{\infty}$ satisfies 
(cf. \cite{L}, \cite{H2}),
\begin{equation*}
(|Rm_{\infty}|^2)_t\le\Delta_{\infty}|Rm_{\infty}|^2
-2|\nabla_{\infty}Rm_{\infty}|^2+16|Rm_{\infty}|^3\quad\mbox{ in }
M_{\infty}\times (0,T], 
\end{equation*} 
by (31),
\begin{align}
(|Rm_{\infty}|^2)_t\le&\Delta_{\infty}|Rm_{\infty}|^2+64|Rm_{\infty}|^2
\quad\mbox{ in }M_{\infty}\times (0,T]\nonumber\\
\Rightarrow\quad
(e^{-64t}|Rm_{\infty}|^2)_t\le&\Delta_{\infty}(e^{-64t}|Rm_{\infty}|^2)
\quad\mbox{ in }M_{\infty}\times (0,T].
\end{align} 
By (35) and the maximum principle \cite{Hs2},
\begin{equation}
e^{-64t}|Rm_{\infty}(x,t)|^2\le e^{-64s_i}\sup_{x\in M}|Rm_{\infty}(x,s_i)|^2
\le h(s_i)^2\quad\mbox{ in }M_{\infty}\times (s_i,T]\quad\forall 
i\in\Z^+. 
\end{equation} 
Letting $i\to\infty$ in (36) we get (34). This contradicts (32).
Hence case (a) cannot occur.

\noindent $\underline{\mbox{\bf Case (b)}}$: $\delta_2>0$.

Let $C_3$ and $C_4$ be as in case 1 and let
\begin{equation}
T_2=\min (T/2,\delta_2/(4C_3)).
\end{equation}
Since $\delta_2>0$, there exists $y_{\infty}\in M$ such that 
\begin{equation}
|Rm_{\infty}|(y_{\infty},T_2)\ge\delta_2/2.
\end{equation} 
Then there exists a sequence $y_k\in M_k=B_{\hat{g}_k(0)}(\2{x}_k,A_k)$ 
such that 
\begin{equation}
d_{\hat{g}_k(T_2)}(y_k,y_{\infty})\to 0\quad\mbox{  as }k\to\infty
\end{equation} 
and 
\begin{equation}
|Rm_k|(y_k,T_2)\to|Rm_{\infty}|(y_{\infty},T_2)\quad\mbox{ as }
k\to\infty.
\end{equation} 
Let $r_0=d_{g_{\infty}(T_2)}(x_{\infty},y_{\infty})+1$. Since 
$d_{\hat{g}_k(T_2)}(\2{x}_k,y_k)\to d_{g_{\infty}(T_2)}(x_{\infty},y_{\infty})$
and $A_k\to\infty$, $Q_k\to\infty$, as $k\to\infty$, by passing to 
a subsequence if necessary we may assume without loss of generality that 
\begin{equation}
d_{\hat{g}_k(T_2)}(\2{x}_k,y_k)<r_0\quad\mbox{ and }\quad
\min (A_k,Q_k^{1/2}/2)>1+e^{3nT}r_0\quad\forall k\in\Z^+.
\end{equation} 
Then by (25) and (41),
\begin{equation}
B_{\hat{g}_k(t)}(\2{x}_k,r_0)\subset B_{\hat{g}_k(0)}(\2{x}_k,e^{3nT}r_0)
\quad\forall 0\le t\le\hat{t}_k,k\in\Z^+.
\end{equation} 
By (41) and (42),
\begin{equation}
B_{\hat{g}_k(0)}(y_k,1)\subset B_{\hat{g}_k(0)}(\2{x}_k,e^{3nT}r_0+1)
\subset B_{\hat{g}_k(0)}(\2{x}_k,\min (A_k,Q_k^{1/2}/2))
\end{equation}
holds for any $k\in\Z^+$.
By (13), (43), and an argument similar to the proof of Case 1 but with
$(y_k,T/2)$ replacing $(\2{x}_k,\hat{t}_k)$ in the proof there we get
\begin{equation}
|Rm_{\hat{g}_k(T_2)}|(y_k,T_2)|\le Q_k^{-2}e^{C_4T_2}+C_3T_2\quad\forall
k\in\Z^+.
\end{equation} 
Letting $k\to\infty$ in (44), by (37), (38) and (40), we get
\begin{equation*}
\delta_2/2\le C_3T_2\le\delta_2/4
\end{equation*} 
and contradiction arises. Hence case (b) is false. Thus no such 
sequence of manifolds $(M_k,g_k,x_{0,k})$ exists and the theorem 
follows. 

\hfill$\square$\vspace{6pt}

By the result of \cite{H3} and an argument similar to the proof 
of case 2 we have the following extension of the compactness 
result of Hamilton \cite{H3}.

\begin{thm}
Let $\{(M_k,g_k(t),x_k)\}_{k=1}^{\infty}$ be a sequence of complete
pointed Riemannian manifolds evolving under Ricci flow (1) 
with sectional curvatures uniform bounded above by some constant 
$B>0$ on $[0,T]$ and uniform positive lower bound on the 
injectivity radii at $x_k$ with respect to the metric $g_k(0)$. 
Then there exists a subsequence which converges to some complete 
pointed Riemannian manifold $(M,g(t),x_0)$ on $(0,T]$ that evolves 
by the Ricci flow on $(0,T]$ with  sectional curvatures uniform 
bounded above by $B$.
\end{thm}


\end{document}